\RequirePackage{fix-cm}
\documentclass[smallcondensed]{svjour3}     
\smartqed  
\usepackage{graphicx}
\begin{document}

\title{Some  Ricci-Flat ($\alpha,\beta$)-Metrics%
\thanks{Authors are both supported in part by The Scientific and Technological Research Council of Turkey (TUBITAK), Grant (No. 113F311)}
}


\author{Esra Sengelen Sevim         \and
        Semail \"Ulgen 
}


\institute{Esra Sengelen Sevim \at
              Department of Mathematical Sciences\\
                \.{I}stanbul Bilgi University\\
            Eski Silahtaraga Elektrik Santrali\\
                Kazim Karabekir Cad. No: 2/13\\
                34060 Ey\"up, Istanbul\\
              Tel.: +90-212-3115425\\
             \email{esra.sengelen@bilgi.edu.tr}           
           \and
            Semail \"Ulgen \at
              Antalya International  University\\
            \"Universite Cad. No:2 07190\\
            D\"osemealti, Antalya\\
             Tel.: +90-242-2450185\\
            \email{sulgen@antalya.edu.tr}
}

\date{Received: date / Accepted: date}

\maketitle

\begin{abstract}
In this paper, we study a special class of Finsler metrics, $(\alpha,\beta)$-metrics, defined by $F=\alpha\phi(\beta/\alpha)$, where $\alpha$ is a Riemannian metric and $\beta$ is a 1-form. We find an equation that characterizes Ricci-flat $(\alpha,\beta)$-metrics under the condition that the length of $\beta$ with respect to $\alpha$ is constant.
\end{abstract}

\section{Introduction}
\label{intro}
Riemannian metrics on a manifold are quadratic metrics, while Finsler metrics are those without restriction on the quadratic property. The Riemannian curvature in Riemannian geometry can be extended to Finsler metrics as a family of linear transformations on the tangent spaces. The Ricci curvature is  the trace of the Riemann curvature. It is a natural problem to study Finsler metrics with isotropic Ricci curvature ${Ric}={Ric}(x,y)$ and
\begin{equation} {Ric}=(n-1)\tau F^2\nonumber\end{equation}
where $\tau=\tau(x)$ is a scalar function on the n-dimensional manifold and $F(x,y)$ is a Finsler metric. Such metrics are called Einstein Finsler metrics.

In this paper, we consider Einstein metrics defined by a Riemannian metric $\alpha$ and 1-form $\beta$ in the following form:
 \begin{equation} F=\alpha\phi(s),~s={\beta\over\alpha},\label{alphabeta}\end{equation}
where $\phi=\phi(s)$ is a positive smooth function. Finsler metrics defined in (\ref{alphabeta}) are called $(\alpha,\beta)$-metrics.

The simplest $(\alpha,\beta)$-metrics are Randers metrics also defined by $F=\alpha+\beta$. In \cite{BaRo}, Bao-Robles find equations on $\alpha$ and $\beta$ that characterize Randers metrics of constant Ricci curvature. There are many Randers metrics of constant Ricci curvature. Thus one just needs to focus on Ricci-flat $(\alpha,\beta)$-metrics. In \cite{CShenT} and \cite{ESZ}, the authors obtained equations on $\alpha$, $\beta$ and $\phi$ that characterize Ricci-flat $(\alpha,\beta)$-metrics of Douglas type. In \cite{ESZ1}, the authors obtained equations on $\alpha$, $\beta$ and $\phi$ that characterize Ricci-flat $(\alpha,\beta)$-metrics which is not of Douglas type. In this paper, we show that there are some more Ricci-flat $(\alpha,\beta)$-metrics.\\

In this paper, we prove the following theorem.\\

\textbf{Theorem 1.1 }\label{thmES1.1}
Let  $F = \alpha \phi(s), s={\beta}/{\alpha}$ be an $(\alpha,\beta)$-metric on an $n$-dimensional manifold $M$ where $\alpha=\sqrt{a_{ij}y^iy^j}$ is a Riemannian metric, $\beta=b_iy^i$ is a 1-form  and $\phi=\phi(s)$ is a positive $C^{\infty}$ function. Suppose that $\alpha$, $\beta$ and $\phi$ satisfy the following conditions:

 \begin{eqnarray}
   (a)&~& ^\alpha \mathbf{Ric}=(n-1)(c_1\alpha^2+c_2\beta^2)\tau,\nonumber\\
   (b)&~& r_{ij}=0,\nonumber\\
   (c)&~&  s_j=0,\nonumber\\
   (d)&~&  t_{ij}=(c_1+c_2b^2)(b_ib_j-a_{ij}b^2)\tau,\nonumber\\
   (e)&~& \phi~~ satisfies \nonumber\\
   &~&  0=(c_1+c_2s^2)+(c_1+c_2b^2)\Big \{2\frac{(s^2-b^2)}{(n-1)}(Q^{\prime}-Q^2+sQQ^{\prime})+Q^2b^2+ 2Qs\Big\},\nonumber\\
&~ &\label{lemma1}\end{eqnarray}
where $b:=\sqrt{a^{ij}b_ib_j}$, $c_1$ and $c_2$ are constants, $\tau:=\tau(x)$ is a scalar function, $t_{ij}:=s_{im}s^m_j$ and
 $$  Q:= \frac{\phi'}{\phi-s\phi'}, $$
Then F is Ricci-flat.
\\
\\
The equation (\ref{lemma1}) is an ordinary differential equation. It is of first order in $Q$ and second order in $\phi$. According to the ODE theory,  the local solution of (\ref{lemma1}) exists nearby $s=0$ for any given initial conditions. But we are unable to express it in terms of elementary functions and we are unable to show that the solution is defined on an interval containing $[-b,b]$ . Thus the $(\alpha,\beta)$-metric $F=\alpha\phi(\beta/\alpha)$ defined by $\phi$ might be singular. We can give the the following example taking $c_2=0$ in Theorem 1.1, then $\alpha$, $\beta$ satisfies Theorem 1.1 (a)-(e). Then for any $\phi=\phi(s)$ satisfying (\ref{lemma1}), we obtain a (possibly singular)  Ricci-flat $(\alpha,\beta)$-metrics.\\

\textbf{Example 1.1.}
Let $F=\alpha+\beta$ be the family of Randers metrics on $S^3$ constructed in \cite{BaRo1}(see also \cite{ZShen}). It is shown that $r_{ij}=0$ and $s_j=0$. Thus for any $C^{\infty}$ positive function $\phi=\phi(s)$ satisfying (\ref{alphabetadefin}), the $(\alpha,\beta)$-metric $F=\alpha\phi(\beta/\alpha)$ has vanishing $S$-curvature.

\section{Preliminaries}
\label{sec:1}
A Finsler metric on a manifold $M$ is a nonnegative scalar function $F=F(x, y)$ on the tangent bundle $TM$, where $x$ is a point in $M$ and $y\in T_xM$ is a tangent vector at $x$.
In local coordinates, the geodesics of a Finsler metric $F=F(x,y)$
are characterized by
$${d^2x^i\over dt^2}+2G^i(x,{dx\over dt})=0,$$
where \begin{equation} G^i:={1\over 4}g^{il}(x,y)\Big \{[F^2]_{x^ky^l}(x,y)y^k-[F^2]_{x^l}(x,y)\Big \}, \label{P(2.1)} \end{equation}
where $g_{ij}=\frac{1}{2}[F^2]_{y^iy^j}$. The local functions $G^i$ on $TM$ define a global vector field
\[ G = y^i \frac{\partial }{\partial x^i} - 2 G^i \frac{\partial }{\partial y^i}.\]
The vector field $G$ is called the {\it spray} of $F$ and the local functions
$G^i=G^i(x,y)$ are called  {\it spray coefficients} of $F$.

For any $x\in M$ and $y\in T_xM \backslash \{0\}$, the Riemann
curvature $\textbf{R}_y: T_xM\to T_xM$ is defined by ${\bf R}_y(u)= R^i_{\ k}(x,y) u^k \frac{\partial}{\partial x^i}|_x$, where
$$ R^i_{\ k}=2{\partial G^i\over\partial
x^k}-{\partial^2 G^i\over\partial x^m\partial
y^k}y^m+2G^m{\partial^2 G^i\over\partial y^m\partial y^k}-{\partial
G^i\over\partial y^m}{\partial G^m\over\partial y^k}.$$
Then the Ricci
curvature is given by
$$\textbf{Ric}=2{\partial G^i\over\partial
x^i}-{\partial^2 G^i\over\partial x^m\partial
y^i}y^m+2G^m{\partial^2 G^i\over\partial y^m\partial y^i}-{\partial
G^i\over\partial y^m}{\partial G^m\over\partial y^i}.$$

An $(\alpha,
\beta)$-metric on a manifold $M$ is a scalar function on $TM$ defined by $$F:=\alpha\phi(s),\ \ \ \ \
s=\frac{\beta}{\alpha},$$ where $\phi=\phi(s)$ is a $C^\infty$ function on
$(-b_0,b_0)$, $\alpha=\sqrt{a_{ij}(x)y^iy^j}$ is a Riemannian metric and
$\beta=b_i(x)y^i$ is a 1-form with $b(x):=\|\beta_x\|_{\alpha} < b_0$. It can be shown that for any Riemannian
metric $\alpha$ and any 1-form $\beta$ on $M$ with $b(x)<b_0$ the
function $F=\alpha\phi(\beta/\alpha)$ is a (positive definite) Finsler metric
if and only if $\phi$ satisfies
\begin{equation} \phi(s)>0,\ \ \ \ \phi(s)-s\phi'(s)+(\rho^2-s^2)\phi''(s)>0,\ \ \ \ (|s|\leq\rho<  b_0).\label{alphabetadefin}\end{equation}
Let
\begin{eqnarray*}
&&r_{ij}:=\frac{1}{2}(b_{i|j}+b_{j|i}),\ \ \ \ \ \ s_{ij}:=\frac{1}{2}(b_{i|j}-b_{j|i}),\\
&&r_j:=b^ir_{ij},\ \ \ \ \ \ \ \ \ \ \ \ \ \ \ \ \ \ s_j:=b^is_{ij},
\end{eqnarray*}
where "$|$" denotes the covariant derivative with respect to the
Levi-Civita connection of $\alpha$.
By (\ref{P(2.1)}), the spray
coefficients $G^i$ of $F$ are given by the following Lemma.
\\
\\

\textbf{Lemma 2.1.} (\cite{ZShen1})
For an $(\alpha,\beta)$-metric $F = \alpha\phi(s), s=\beta/\alpha$, the spray
coefficients of $F$ are given by \begin{equation} G^i={^\alpha G^i}+\alpha Q s^i_{\
0}+\Theta \{r_{00}-2Q\alpha s_0\}\frac{y^i}{\alpha}+\Psi\{r_{00}-2Q\alpha
s_0\}b^i,\label{L(2.1)} \end{equation} where $^\alpha G^i$ are the spray coefficients of $\alpha$,
\begin{eqnarray*}
Q&:=&\frac{\phi'}{\phi-s\phi'},\\
\Theta&:=& \frac{Q-sQ'}{2\Delta},\\
\Psi&:=&\frac{Q'}{2\Delta},\\
\Delta & : = & 1+sQ+(b^2-s^2)Q'
\end{eqnarray*}
and $s^i_{\ j}:=a^{ik}s_{kj}$, $s_{ij}:=a_{ih}s^h_j$. The index "0" means contracting with
$y$, for example, $s^i_{\ 0}:=s^i_{\ j}y^j, s_0:=s_iy^i, s_{ij}y^j:=s_{i0}, s_{ij}y^i:=s_{0i},
r_{00}:=r_{ij}y^iy^j. $

\section{Proof of Theorem 1.1}
\label{sec:2}
In this section we prove Theorem 1.1. Throughout this section, we assume that the dimension is greater than two.
First we  give the following Lemma.\\

\textbf{Lemma 3.1. }\label{lemES2.1}
Let $F=\alpha\phi(\beta / \alpha)$ be an $(\alpha,\beta)$-metric on an $n$-dimensional manifold $M$, $n\geq 3$. Suppose that $\alpha=\sqrt{a_{ij}(x)y^iy^j}$ and $\beta=b_i(x)y^i$ satisfy the conditions of Theorem 1.1 (b) and (c), then the following equations are satisfied:
 \begin{eqnarray}
  s^m_{0  |m}&=&(n-1)(c_1+c_2b^2)\tau\beta\label{sm0m}
   \end{eqnarray}
where $\tau=\tau(x)$ is a scalar function and $c_1$ and $c_2$ are constants.
\\

\textbf{Proof:} By Ricci identities, we have

\begin{eqnarray*}
b_{i  | j  | k }-b_{i  | k  | j }&=& b^m{^{^{\alpha}}{\texttt{R}}_{imjk}} ,\\
-b_{k  | i  | j }+b_{k  | j | i }&=&-b^m{^{^{\alpha}}{\texttt{R}}_{kmij}} ,\\
b_{j  | k  | i }-b_{j  | i  | k }&=& b^m{^{^{\alpha}}{\texttt{R}}_{jmki}} .
\end{eqnarray*}
\noindent On the other hand,
\begin{eqnarray*}
b_{i |k  | j }+b_{k  | i  | j }&=& 2r_{ik  |j},\\
-b_{k  | j  | i }-b_{j  | k  | i }&=& -2r_{kj  |i}.
\end{eqnarray*}
\noindent Adding all the equations above, we get
\begin{equation} s_{ij |k}={1\over2}(b_{i |j  | k }+b_{j  | i  | k })=-b^m{^{^{\alpha}}{\texttt{R}}_{kmij}}+r_{ik  |j}-r_{kj  |i}.\nonumber\end{equation}

\noindent The condition $(b)$ in Theorem (\ref{lemma1}) helps one to  rewrite the above equation as follows:
\begin{equation} s_{ij |k}={1\over2}(b_{i |j  | k }+b_{j  | i  | k })=-b^m{^{^{\alpha}}{\texttt{R}}_{kmij}}.\label{sm0m1}\end{equation}
Hence,
\begin{equation} s^m_{0  | m}=b^m{^{^{\alpha}}{\texttt{Ric}}}_{m0}+r^m_{m|0}-r^m_{0|m}.\label{sm0m00}\end{equation}
The condition (b) in Theorem 1.1, (\ref{sm0m00}) implies the following:
\begin{equation} s^m_{0  | m}=b^m{^{^{\alpha}}{\texttt{Ric}}}_{m0}.\label{sm0m2}\end{equation}

\noindent The condition $(a)$ in Theorem (\ref{lemma1}) implies the following
\begin{eqnarray}
{^{^{\alpha}}}{\texttt{Ric}}_{ij}&=&\Big( {1\over 2} {^{^{\alpha}}}\mathbf{Ric}\Big)_{y^iy^j}\nonumber\\
&=& (n-1)\Big( c_1a_{ij}+c_2b_ib_j\Big)\tau,\label{sm0m3}
\end{eqnarray}
\noindent and we obtain:
\begin{equation} b^m{^{^{\alpha}}{\texttt{Ric}}}_{m0}=(n-1)(c_1+c_2b^2)\tau\beta.\label{sm0m4}\end{equation}

\noindent Hence, the equation in (\ref{sm0m}) follows from equations (\ref{sm0m2}) and (\ref{sm0m4}).

\qed
\vskip 0,2cm

Next, we compute the Ricci curvature of the $(\alpha,\beta)$-metric under the conditions $(a)-(e)$ of the Theorem 1.1.
By Lemma 2.1, the spray coefficients of $F$ can be written as
\begin{equation}\nonumber G^i={~}^{\alpha}G+T^i,\nonumber\end{equation}
where
\begin{equation} T^i=\alpha Q s^i_0.\nonumber\end{equation}

It is well known (\cite{ZShen1}) that the curvature tensor can be written as
\begin{equation} R^i_k=^\alpha R^i_k+H^i_k,\nonumber\end{equation}
where
\begin{equation} H^i_k:=2T^i_{ |k}-T^i_{ |j\cdot k}y^j+2T^jT^i_{\cdot j\cdot k}-T^i_{\cdot j}T^j_{\cdot k},\nonumber\end{equation}
and $"."$ and $"|"$ mean vertical covariant derivative and horizontal covariant derivative with respect to $\alpha$, respectively. Then
\begin{equation}{\mathbf{Ric}}={~}^{\alpha}{\mathbf{Ric}}+H^i_i,\nonumber\end{equation}
where ${~}^{\alpha}{{\mathbf{Ric}}}$ denotes the Ricci curvature of $\alpha$ and
\begin{equation} H^i_i:=2T^i_{ |i}-T^i_{ |j\cdot i}y^j+2T^jT^i_{\cdot j\cdot i}-T^i_{\cdot j}T^j_{\cdot i}.\label{H}\end{equation}
\\
\\
To compute the Ricci curvature under the conditions $r_{ij}=0$ and $s_j=0$, we need:
\begin{eqnarray}
&& b_{i |j}=s_{ij},\ \ y_is^i_0=0,\ \ y_is^i_{0 |j}=0, \ \ s_{ij}y^iy^j=0,\nonumber\\
&& y_is^i_{0 |j}=0,\ \ b_is^i_0=0,\ \ b_is^i_j=0,\ \ \ b_is^i_{0 |j}=-s_{ij}s^i_0.\label{defins2}
\end{eqnarray}
We can also easily get
\begin{eqnarray}
&& s_{ \cdot i}={b_i\over\alpha}-s{y_i\over\alpha^2},\ \  s_{ \cdot i}b^i={1\over\alpha}(b^2-s^2), \ \ \nonumber\\
&& s_{ \cdot i}y^i=0,\ \ s_{ \cdot i}s^i_0=0, \ \ s_{\cdot i} s^i_{0 |j}=-\frac{s_{ij}s^i_0}{\alpha}.\label{defins3}
\end{eqnarray}
\begin{eqnarray}
&& s_{\cdot j\cdot i}=-\frac{b_jy_i}{\alpha^3}-\frac{b_iy_j}{\alpha^3}+3s\frac{y_iy_j}{\alpha^4}-s\frac{a_{ji}}{\alpha^2}\nonumber\\
&& s_{\cdot j\cdot i}s^i_0=-\frac{s}{\alpha^2}s_{j0},\ \ s_{\cdot j\cdot i}s^i_0s^j_0=-\frac{s}{\alpha^2}s_{j0}s^j_0.\label{defins}
\end{eqnarray}

\begin{eqnarray}
&&s_{ |i}={s_{0i}\over\alpha},\ \  s_{ |i}y^i=0 \ \  s_{ |i}b^i=0,\ \  s_{|j\cdot i}= \frac{s_{ij}}{\alpha}-\frac{s_{0j}y_i}{\alpha^3}\nonumber\\
&& s_{|j\cdot i}s^i_{0}=\frac{s_{ij}s^i_0}{\alpha},\ \ s_{\cdot j\cdot i}s^i_0s^j_0=-\frac{s}{\alpha^2}s_{j0}s^j_0.\label{defins1}
\end{eqnarray}

\noindent Using the above identities in   (\ref{defins1}),   the equation   $T^i_{|i}=\alpha Q^{\prime}s_{|i} s^i_0+\alpha Q s^i_{0|i}$       is simplified to


\begin{eqnarray}
T^i_{ |i}&=&Q^{\prime}s_{0i}s^i_{0}+\alpha Q s^i_{0 |i}.\label{eq1}
\end{eqnarray}

\noindent The identities in  (\ref{defins}) and (\ref{defins1}) are used in
\begin{eqnarray*}
T^i_{|j}&=&\alpha Q^{\prime}s_{|j} s^i_0+\alpha Q s^i_{0|j},\\
T^i_{|j\cdot i}&=&\frac{y_i}{\alpha} Q^{\prime}s_{|j} s^i_0+\alpha Q^{\prime\prime}s_{\cdot i}s_{|j} s^i_0+\alpha Q^{\prime}s_{|j\cdot i} s^i_0+\alpha Q^{\prime}s_{|j} s^i_i+\frac{y_i}{\alpha} Q s^i_{0 |j}+\alpha Q^{\prime}s_{\cdot i} s^i_{0 |j}+\alpha Q s^i_{i|j},
\end{eqnarray*}
to get  the following simplified  equations:
\begin{eqnarray}
T^i_{|j\cdot i}&=&\alpha Q^{\prime}\frac{s_{ij}s^i_0}{\alpha}
-\alpha Q^{\prime}\frac{s_{ij}s^i_0}{\alpha},\nonumber\\
T^i_{|j\cdot i}&=&0,\nonumber\\
T^i_{|j\cdot i}y^j&=&0.\label{eq2}\end{eqnarray}

\noindent We further have
\begin{eqnarray*}
T^i_{\cdot j}&=&\frac{y_j}{\alpha} Q s^i_0+\alpha Q^{\prime}s_{\cdot j} s^i_0+\alpha Q s^i_j,\\
T^i_{\cdot j\cdot i}&=&\Big(\frac{a_{ij}}{\alpha}-\frac{y_iy_j}{\alpha^3}\Big)Q s^i_0+\frac{y_j}{\alpha} Q^{\prime}s_{\cdot i} s^i_0+\frac{y_j}{\alpha} Q s^i_i+\frac{y_i}{\alpha} Q^{\prime}s_{\cdot j} s^i_0+\alpha Q^{\prime\prime}s_{\cdot i}s_{\cdot j} s^i_0\\
&&+\alpha Q^{\prime}s_{\cdot j\cdot i} s^i_0+\alpha Q^{\prime}s_{\cdot j} s^i_i
+\frac{y_i}{\alpha} Q s^i_j+\alpha Q^{\prime}s_{\cdot i} s^i_j,\\
T^jT^i_{\cdot j\cdot i}&=&\alpha Q s^j_0\Big\{\Big(\frac{a_{ij}}{\alpha}-\frac{y_iy_j}{\alpha^3}\Big)Q s^i_0+\frac{y_j}{\alpha} Q^{\prime}s_{\cdot i} s^i_0+\frac{y_j}{\alpha} Q s^i_i+\frac{y_i}{\alpha} Q^{\prime}s_{\cdot j} s^i_0+\alpha Q^{\prime\prime}s_{\cdot i}s_{\cdot j} s^i_0\\
&&+\alpha Q^{\prime}s_{\cdot j\cdot i} s^i_0+\alpha Q^{\prime}s_{\cdot j} s^i_i+\frac{y_i}{\alpha} Q s^i_j+\alpha Q^{\prime}s_{\cdot i} s^i_j\Big\}.
\end{eqnarray*}
Using the identities in (\ref{defins2}), (\ref{defins3}), (\ref{defins}) and (\ref{defins1}), we get:
\begin{eqnarray*}
T^jT^i_{\cdot j\cdot i}&=&Q^2s_{i0}s^i_0-sQQ^{\prime}s_{j0}s^j_{0}+Q^2s^j_0s^0_j-sQQ^{\prime}s^j_{0}s^0_{j}.
\end{eqnarray*}
Using the fact that $s^0_j=-s_{j0}$, we obtain the following simple equation:
\begin{eqnarray}
T^jT^i_{\cdot j\cdot i}&=&0\label{eq3}.
\end{eqnarray}
After multiplying the following equations:
\begin{eqnarray*}
T^i_{\cdot j}&=&\frac{y_j}{\alpha}Qs^i_0+\alpha Q^{\prime}s_{\cdot j}s^i_0+\alpha Qs^i_j,\\
T^j_{\cdot i}&=&\frac{y_i}{\alpha}Qs^j_0+\alpha Q^{\prime}s_{\cdot i}s^j_0+\alpha Qs^j_i,
\end{eqnarray*}
and then simplifying them we get
\begin{eqnarray}
T^i_{\cdot j}T^j_{\cdot i}&=&2Q^2s_{0i}s^i_0-2sQQ^{\prime}s_{0i}s^i_0+\alpha^2Q^2s^i_js^j_i.\label{eq4}
\end{eqnarray}

\noindent Plugging (\ref{eq1}), (\ref{eq2}), (\ref{eq3}) and (\ref{eq4}) into (\ref{H}), we obtain
\begin{equation} H^i_i=2(Q^{\prime}-Q^2+sQQ^{\prime})t_{00}-{\alpha}^2Q^2t^m_m+2\alpha Q s^i_{0 |i},\label{eqH}\end{equation}
where $t_{00}=t_{ij}y^iy^j$, $t_{00}=(c_1+c_2b^2)(s^2-b^2)\tau\alpha^2.$  Hence $H^i_i$ and also $ \mathbf{Ric}$ are expressed as follows:
\begin{equation} H^i_i=(c_1+c_2b^2)\tau\Big \{2(Q^{\prime}-Q^2+sQQ^{\prime})(s^2-b^2)+(n-1)Q^2b^2+ 2(n-1)sQ\Big\}\alpha^2. \nonumber\end{equation}
and
\begin{equation} \mathbf{Ric}={~}^{\alpha}{\mathbf{Ric}}+\tau\alpha^2\Gamma\nonumber \end{equation}
where
\begin{equation}\Gamma:= (c_1+c_2b^2)\Big \{2(Q^{\prime}-Q^2+sQQ^{\prime})(s^2-b^2)+(n-1)Q^2b^2+ 2(n-1)sQ\Big\}.\nonumber\end{equation}
Thus $\mathbf{Ric}=0$ if and only if
\begin{equation} (n-1)(c_1+c_2s^2)\tau\alpha^2+\Gamma=0\label{eq5}\end{equation}
We can rewrite (\ref{eq5}) as (\ref{lemma1}).
\qed

\end{document}